\documentclass[11pt]{article}
\usepackage{amsmath,amssymb,a4wide}

\def\Ksep{{K^{\rm sep}}}

 \def\Aut{\mathop{\rm Aut}\nolimits}

 \def\Gal{\mathop{\rm Gal}\nolimits}
 \def\End{\mathop{\rm End}\nolimits}

\def\Stab{\mathop{\rm Stab}\nolimits}
\def\Fix{\mathop{\rm Fix}\nolimits}

\def\codim{\mathop{\rm codim}\nolimits}

\def\rank{\mathop{\rm rank}\nolimits}

\def\GL{\mathop{\rm GL}\nolimits}

\def\GSp{\mathop{\rm GSp}\nolimits}

\def\tor{{\rm tor}}
\def\sep{{\rm sep}}

\let\phi\varphi
\let\epsilon\varepsilon

\newtheorem{Thm}{Theorem}[section]
\newtheorem{Prop}[Thm]{Proposition}
\newtheorem{Lem}[Thm]{Lemma}

\newtheorem{Cor}[Thm]{Corollary}

\def\qed{{\hskip0pt\unskip\unskip\nobreak\hfil\penalty50
          \hskip1em\hbox{}\nobreak\hfil
           {$\square$}
          \parfillskip=0pt\finalhyphendemerits=0
          \par}\medskip}

\newenvironment{Proof}
               {\noindent{\bf Proof.}\ }
               {\qed}
               
               {\noindent{\bf Proof sketch.}\ }
               {\qed}               

\newenvironment{Proofof}[1]
               {\noindent{\bf Proof of #1.}\ }
               {\qed}

\newcommand{\BG}{{\mathbb{G}}}

\newcommand{\BN}{{\mathbb{N}}}

\newcommand{\BQ}{{\mathbb{Q}}}

\newcommand{\BZ}{{\mathbb{Z}}}

\newcommand{\Fa}{{\mathfrak{a}}}
\newcommand{\Fb}{{\mathfrak{b}}}

\newcommand{\Fp}{{\mathfrak{p}}}

\def\Gm{\BG_{\mathrm m}}

\newbox\mybox
\def\arrover#1{\mathrel{
       \setbox\mybox=\hbox spread 1.4em
              {\hfil$\scriptstyle#1$\hfil}
       \vbox{\offinterlineskip\copy\mybox
             \hbox to\wd\mybox{\rightarrowfill}}}}

\def\larrover#1{\mathrel{
       \setbox\mybox=\hbox spread 1.4em
              {\hfil$\scriptstyle#1\vphantom{g}$\hfil}
       \vbox{\offinterlineskip\copy\mybox
             \hbox to\wd\mybox{\leftarrowfill}}}}

\def\ontoover#1{\mathrel{
       \setbox\mybox=\hbox spread 1.4em
              {\hfil$\scriptstyle#1\vphantom{g}$\hfil}
       \vbox{\offinterlineskip\copy\mybox
             \hbox to\wd\mybox{\rightarrowfill\hskip-2.8mm
                               $\rightarrow$}}}}
\def\leftontoover#1{\mathrel{
       \setbox\mybox=\hbox spread 1.4em
              {\hfil$\scriptstyle#1\vphantom{g}$\hfil}
       \vbox{\offinterlineskip\copy\mybox
             \hbox to\wd\mybox{$\leftarrow$\hskip-2.8mm
                               \leftarrowfill}}}}
\let\longto\longrightarrow

\begin{document}

\title{Torsion bounds for elliptic curves and Drinfeld modules}
\author{Florian Breuer (fbreuer@sun.ac.za) \\ Stellenbosch University, Stellenbosch, South Africa}
\maketitle

\begin{abstract}
We derive asymptotically optimal upper bounds on the number of $L$-rational torsion points on a given elliptic curve or Drinfeld module defined over a finitely generated field $K$, as a function of the degree $[L:K]$. Our main tool is the adelic openness of the image of Galois representations attached to elliptic curves and Drinfeld modules, due to Serre and Pink-R\"utsche, respectively. Our approach is to prove a general result for certain Galois modules, which applies simultaneously to elliptic curves and to Drinfeld modules.


\end{abstract}

\section{Statement of main result}

Let $A$ be a Dedekind domain whose fraction field $k$ is a global field, and let $M$ be an $A$-module. 
For a non-zero ideal $\Fa\subset A$ we denote by 
\[
M[\Fa]:=\{x\in M\;|\; a\cdot x = 0,\;\forall a\in\Fa\}
\]
the $\Fa$-torsion submodule of $M$, by $M[\Fa^\infty] := \cup_{n\geq 1}M[\Fa^n]$ the $\Fa$-power torsion submodule, and 
by $M_\tor=\cup_{\Fa\subset A} M[\Fa]$ the full torsion submodule.

Let $K$ be a finitely generated field, and denote by $\Ksep$ the separable closure of $K$ in an algebraic closure $\bar{K}$. Let $G_K = \Gal(\Ksep/K)$ act on $M$. For a submodule $H\subset M$ and a finite extension $L/K$ we denote by $H(L)$ the subset of $H$ fixed by $\Gal(\Ksep/L)$. The cardinality of a finite set $S$ is denoted $|S|$.  

The goal of this article is to prove the following.

\begin{Thm}\label{main}
Suppose we are in one of the following two situations:
\begin{itemize}
\item[(a)] $A=\BZ$, $K$ is a finitely generated field of characteristic 0 and $M$ is an elliptic curve over $K$. Set $\gamma=\rank_{\BZ}(\End_{\bar{K}}(M))/2$.
\item[(b)] $k$ is a global function field and $A$ is the ring of elements of $k$ regular outside a fixed place $\infty$ of $k$. 
$M$ is a rank $r$ Drinfeld $A$-module in generic characteristic over the finitely generated field $K$, and set $\gamma=\rank_{A}(\End_{\bar{K}}(M))/r$.
\end{itemize}
Then we have
\begin{itemize}
\item[(I)] Let $\Fp\subset A$ be a non-zero prime ideal. Then there exists a constant $C$ depending on $M,K$ and $\Fp$ such that, for any finite extension $L/K$,
\begin{equation}
\big|M[\Fp^\infty](L)\big| \leq C[L:K]^\gamma.
\end{equation}
\item[(II)] There exists a constant $C$ depending on $M$ and $K$ such that, for any finite extension $L/K$,
\begin{equation}
\big|M_\tor(L)\big| \leq C\big([L:K]\log\log[L:K]\big)^\gamma.
\end{equation}
\end{itemize}
Moreover, these bounds are asymptotically optimal in the sense that there exist towers of fields achieving these bounds for suitable values of $C$.
\end{Thm}



\section{Galois modules}

Continuing with the notation from the previous section, we call $M$ a {\em $(G_K,A)$-module of rank $r$} if the action of $G_K$ commutes with the action of $A$, and $M[\Fa]\cong (A/\Fa)^r$ as $A$-modules for every non-zero ideal $\Fa\subset A$. 

Then for every non-zero ideal $\Fa\subset A$ the action of $G_K$ induces a Galois representation
\[
\rho_{\Fa} : G_K \longto \Aut(M[\Fa]) \cong \GL_r(A/\Fa),
\]
once we have chosen a basis for $M[\Fa]$. We denote the index of the image by
\[
I(\Fa) := \big(\GL_r(A/\Fa) : \rho_{\Fa}(G_K)\big).
\]

Our main tool is

\begin{Thm}\label{GaloisModules}
Let $M$ be a $(G_K,A)$-module of rank $r$.
\begin{itemize}
	\item[(I)] Let $\Fp\subset A$ be a non-zero prime ideal. Suppose that there exists a constant $C_{\Fp}$, depending on $M$, $K$ and $\Fp$, such that $I(\Fp^n) \leq C_{\Fp}$ for all $n\in\BN$. Then there exists a constant $C$ depending on $C_{\Fp}$, $r$ and $K$ such that, for every finite extension $L/K$,
	\[
  \big|M[\Fp^\infty](L)\big| \leq C[L:K]^{1/r}.
	\]  
	\item[(II)] Suppose that there exists a constant $C_0$, depending on $M$ and $K$, such that $I(\Fa)\leq C_0$ for all non-zero $\Fa\subset A$. Then there exists a constant $C$ depending on $C_0$, $r$ and $K$ such that, for every finite extension $L/K$,
	\[
	\big|M_\tor(L)\big| \leq C\big([L:K]\log\log[L:K]\big)^{1/r}.
	\] 
\end{itemize}
Moreover, these bounds are asymptotically optimal in the sense that there exist towers of fields achieving these bounds for suitable values of $C$.
\end{Thm}

\subsection{Elementary Lemmas}

We collect the following elementary results, which we will need in the proof of Theorem \ref{GaloisModules}.
For a non-zero ideal $\Fa\subset A$ we write $|\Fa|:=|A/\Fa|$.
We define the function
\[
\theta(\Fa) := \prod_{\Fp|\Fa}\left(1-\frac{1}{|\Fp|}\right)^{-1},
\]
where the product ranges over all prime ideals $\Fp|\Fa$.

\begin{Lem}\label{log}
There exist constants $C_1,C_2>0$, depending on $A$, such that $\theta(\Fa)\leq C_1\log\log|\Fa|$ for all non-zero $\Fa\subset A$. 
Moreover, if $\Fa_n := \prod_{|\Fp| \leq n}\Fp$, then $\theta(\Fa_n) \geq C_2\log\log|\Fa_n|$ for all $n\in\BN$. 
\end{Lem}

\begin{Proof}
It is clear that $\theta(\Fa)$ achieves its fastest growth (relative to $|\Fa|$) for $\Fa_n=\prod_{|\Fp|\leq n}\Fp$. 

We start with the following version of Mertens' Theorem \cite[Theorems 2 and 3]{Rosen}:
\[
\theta(\Fa_n) = \prod_{|\Fp|\leq n} \left(1-\frac{1}{|\Fp|}\right)^{-1} = C_1\log n + O(1),
\]
for an explicit constant $C_1>0$.
On the other hand, we have,
\[
\log|\Fa_n| = \sum_{|\Fp|\leq n}\log|\Fp| = n + o(1).
\]
When $k$ is a number field, this is \cite[Theorem 2.2]{Rosen}. When $k$ is a function field over the finite field of $q$ elements,
this follows readily from the well-known estimate 
\[
\big|\{\text{$\Fp$ prime} \;|\; |\Fp| = q^m\}\big| = \frac{q^m}{m} + O\left(\frac{q^{m/2}}{m}\right).
\]

\end{Proof}

\begin{Lem}\label{OrderGL}
Let $\Fa\subset A$ be a non-zero ideal. Then 
\[
|\GL_r(A/\Fa)| = |\Fa|^{r^2}\prod_{\Fp|\Fa}\left(1-\frac{1}{|\Fp|}\right)\left(1-\frac{1}{|\Fp|^2}\right)\cdots\left(1-\frac{1}{|\Fp|^r}\right).
\]
\end{Lem}

\begin{Proof}
Since $|\GL_r(A/\Fa)|$ is multiplicative in $\Fa$, it suffices to prove the result for $\Fa=\Fp^e$, where $\Fp\subset A$ is prime.

It is well-known that $|\GL_r(A/\Fp)| = (|\Fp|^r-1)(|\Fp|^r-|\Fp|)\cdots(|\Fp|^r-|\Fp|^{r-1})$, and the general result follows from the exact sequence 
%
\[
1 \to 1 + M_r(\Fp/\Fp^e) \longto \GL_r(A/\Fp^e) \longto \GL_r(A/\Fp) \to 1,
\]
where $M_r$ denotes the additive group of $r\times r$ matrices. 
\end{Proof}

\begin{Lem}\label{composita}
Let $K_i/K$ and $L_i/K_i$ be finite extensions inside $\bar{K}$, for $i=1,2,\ldots,r$. We denote by $\prod_{i=1}^rK_i$ the compositum of the fields $K_1,\ldots, K_r$ inside $\bar{K}$, and similarly for $\prod_{i=1}^rL_i$. Then
\[
\frac{\prod_{i=1}^r[K_i:K]}{[\prod_{i=1}^rK_i :K]} \leq \frac{\prod_{i=1}^r[L_i:K]}{[\prod_{i=1}^rL_i :K]}.
\]
\end{Lem}

\begin{Proof} Elementary.
\end{Proof}

\subsection{Fields of definition}

Let $H\subset M[\Fa]$ be a subset, define 
\[
\Fix_{\Aut(M[\Fa])}(H) := \{\sigma\in\Aut(M[\Fa]) \;|\; \sigma(h)=h,\;\forall h\in H\},
\] 
and denote by $K(H)$ the field generated by $H$ over $K$, i.e. $K(H)$ is the fixed field of 
$\rho_{\Fa}^{-1}\big(\Fix_{\Aut(M[\Fa])}(H)\big)$. When $H=\{x\}$ we write $K(H)=K(x)$. 
Then $\rho_{\Fa}$ induces an isomorphism:
\begin{Lem}\label{Stab} 
\[
\Gal\big(K(H)/K\big) \cong \frac{\rho_{\Fa}(G_K)}{\rho_{\Fa}(G_K)\cap\Fix_{\Aut(M[\Fa])}(H)}.
\]
\qed
\end{Lem}

Let $x\in M_{\tor}$. Then we say that $x$ has order $\Fa$, where $\Fa\subset A$ is a non-zero ideal, if $x\in M[\Fa]$ but $x\not\in M[\Fb]$ for any ideal $\Fb\supsetneq \Fa$. 
We also denote by $\zeta_A(s) = \prod_{\Fp}\big(1-|\Fp|^{-s}\big)^{-1}$ the zeta-function of $A$.

\begin{Prop}\label{Kx}
Let $x\in M_\tor$ be a point of order $\Fa$.
Then 
\[
[K(x):K] = \frac{1}{C}|\Fa|^r\prod_{\Fp|\Fa}\left(1-\frac{1}{|\Fp|^r}\right)
\]
where $1\leq C \leq I(\Fa)$.
\end{Prop}

\begin{Proof}
Choose a basis for $M[\Fa]$ such that $x$ is the first basis element. This choice determines the isomorphism $\Aut(M[\Fa])\cong\GL_r(A/\Fa)$. The stabilizer of $x$ in $\GL_r(A/\Fa)$ is of the form
\[
\Fix_{\GL_r(A/\Fa)}(x) = \left(\begin{array}{l|l} 1 & * \quad\cdots\quad * \\ \hline 0 & \\ \vdots & \GL_{r-1}(A/\Fa) \\ 0 &    
\end{array}\right)
\]
where the starred entries of the first row are arbitrary elements of $A/\Fa$ and the bottom right $(r-1)\times(r-1)$ block is $\GL_{r-1}(A/\Fa)$. It follows from Lemma \ref{OrderGL} that
\begin{eqnarray*}
|\Fix_{\GL_r(A/\Fa)}(x)| & = & |\Fa|^{r-1}|\GL_{r-1}(A/\Fa)| \\ 
& = & |\Fa|^{r(r-1)}\prod_{\Fp|\Fa}\left(1-\frac{1}{|\Fp|}\right)\left(1-\frac{1}{|\Fp|^2}\right)\cdots\left(1-\frac{1}{|\Fp|^{r-1}}\right).
\end{eqnarray*}

From Lemma \ref{Stab} follows that
\begin{eqnarray*}
[K(x):K] 
& = & \frac{1}{C} \frac{|\GL_r(A/\Fa)|}{|\Fix_{\GL_r(A/\Fa)}(x)|} \\
& = & \frac{1}{C} |\Fa|^r \prod_{\Fp|\Fa}\left(1-\frac{1}{|\Fp|^r}\right),
\end{eqnarray*}
where $1\leq C \leq I(\Fa)$.
\end{Proof}

\begin{Prop}\label{Indep}
Suppose $r\geq 2$.
Then there exists a constant $C>0$ depending on $I(\Fa)$ and $K$, such that the following holds.
Let $x_1,\ldots,x_r\in M[\Fa]$ be a basis for $M[\Fa]$. Then
\[
\frac{\prod_{i=1}^r[K(x_i):K]}{[K(M[\Fa]):K]} \leq C\,\theta(\Fa). 
\] 
\end{Prop}

Of course, $[K(x_1):K] = [K(M[\Fa]):K]$ if $r=1$.

\medskip

\begin{Proof}
From Lemma \ref{Stab} we obtain
\[
[K(M[\Fa]):K] = \frac{1}{I(\Fa)}|\GL_r(A/\Fa)|.
\]
Now Lemma \ref{OrderGL} and Proposition \ref{Kx} give
\begin{eqnarray*}
\frac{\prod_{i=1}^r[K(x_i):K]}{[K(M[\Fa]):K]}  & \leq & I(\Fa)
\prod_{\Fp|\Fa}\left(1-\frac{1}{|\Fp|^r}\right)^{r-1}\left(1-\frac{1}{|\Fp|}\right)^{-1}\left(1-\frac{1}{|\Fp|^2}\right)^{-1}\cdots\left(1-\frac{1}{|\Fp|^{r-1}}\right)^{-1}  \\ 
& \leq &  I(\Fa)\,\zeta_A(2)\zeta_A(3)\ldots\zeta_A(r-1) \cdot\prod_{\Fp|\Fa}\left(1-\frac{1}{|\Fp|}\right)^{-1} 
\end{eqnarray*}
\end{Proof}

The intuition is that the fields generated by linearly independent torsion points have minimal intersection. Explicitly,

\begin{Cor}
There exists a constant $C>0$ depending on $I(\Fa)$ and $K$, such that the following holds.
Let $x_1,x_2\in M[\Fa]$ be points of order $\Fa$ for which $\langle x_1\rangle \cap \langle x_2 \rangle = \{0\}$. Then
\[
[K(x_1)\cap K(x_2) : K] \leq C \theta(\Fa).
\] 
\qed
\end{Cor}

\subsection{Proof of Theorem \ref{GaloisModules}}

\begin{Proofof}{Theorem \ref{GaloisModules}}
Let $H = M_\tor(L)$, which is finite since we assume that the indices $I(\Fa)$ are bounded. Let $\Fa\subset A$ be the minimal ideal for which $H\subset M[\Fa]$. One can choose a basis $x_1,\ldots,x_r$ of $M[\Fa]\cong (A/\Fa)^r$ such that $H = \langle y_1,\ldots,y_r\rangle$, with $y_i \in\langle x_i\rangle$, and $y_i$ is of order $\Fa_i$, for each $i=1,\ldots,r$. Then $K(H)$ is the compositum of the $K(y_i)$'s in $L$.

From Lemma \ref{composita}, Proposition \ref{Indep} and Lemma \ref{log} we obtain 
\[
\frac{\prod_{i=1}^r[K(y_i):K]}{[K(H):K]} \leq \frac{\prod_{i=1}^r[K(x_i):K]}{[K(M[\Fa]):K]} \leq C_1\theta(\Fa) \leq C_2\log\log|\Fa|,
\]
for some constant $C_2$ independent of $H$. From Proposition \ref{Kx} now follows that 
\begin{eqnarray*}
[K(H):K] & \geq & \frac{\prod_{i=1}^r[K(y_i):K]}{C_2\log\log|\Fa|} \qquad\qquad\quad\;\text{(or}\quad
[K(y_1):K] \quad\text{if $r=1$)} \\
& \geq & \frac{1}{I(\Fa)^r\zeta_A(r)^r}\frac{\prod_{i=1}^r |\Fa_i|^r}{C_2\log\log|\Fa|} \qquad\text{(or}\quad \frac{1}{I(\Fa)}\frac{|\Fa_1|}{C_2\log\log|\Fa|}\quad\text{if $r=1$)} \\
& \geq & C_2\frac{|H|^r}{\log\log|H|} ,
\end{eqnarray*}
where $C_2$ is independent of $H$, by the assumption on $I(\Fa)$.

It follows that $|H| \leq C_3\big([K(H):K]\log\log[K(H):K]\big)^{1/r}$, which proves part (II).

If $H=M[\Fp^\infty](L)$, then in the above argument we find that $\Fa=\Fp^n$ for some $n$, and $\theta(\Fp^n)=(1-|\Fp|^{-1})^{-1}$ only depends on $\Fp$, so the $\log\log$-term falls away. Part (I) follows. 
%

Lastly, we show that the bounds are sharp. In case (II), let $\Fa_n := \prod_{|\Fp|\leq n}\Fp$ for $n\in\BN$, as in  Lemma \ref{log}. Now set $L_n:=K(M[\Fa_n])$. By Lemmas \ref{Stab} and \ref{OrderGL},
\begin{eqnarray*}
[L_n:K] & = & \frac{1}{I(\Fa_n)}|\GL_r(A/\Fa_n)| = \frac{1}{I(\Fa_n)}|\Fa_n|^{r^2}\prod_{\Fp|\Fa_n}
\left(1-\frac{1}{|\Fp|}\right)\left(1-\frac{1}{|\Fp|^2}\right)\cdots\left(1-\frac{1}{|\Fp|^r}\right). \\
& \leq & \frac{|\Fa_n|^{r^2}}{\theta(\Fa_n)} \leq C_1\frac{|\Fa_n|^{r^2}}{\log\log(|\Fa_n|)}.
\end{eqnarray*}
It follows that 
\[
|M_\tor(L_n)| \geq \big|M[\Fa_n]\big| = |\Fa_n|^r \geq C_2 \big([L_n:K]\log\log[L_n:K]\big)^{1/r}.
\]
Case (I) is similar: We let $\Fa_n := \Fp^n$ and $L_n:=K(M[\Fp^n])$. This time $\theta(\Fp^n)$ is constant and we find
\[
|M[\Fp^\infty](L_n)| \geq \big|M[\Fp^n]\big| = |\Fp^n|^r \geq C_2 [L_n:K]^{1/r}.
\]
\end{Proofof}

\section{Proof of the main result}

\subsection{Drinfeld modules}

Suppose that $k$ is a global function field, and fix a place $\infty$ of $k$. Let $A$ be the ring of elements of $k$ regular away from $\infty$, and let $\phi$ be a Drinfeld $A$-module of rank $r$ in generic characteristic defined over the finitely generated field $K$. We denote by $\End_L(\phi)$ the ring of endomorphisms of $\phi$ defined over a field $L/K$. See \cite[Chapter 4]{Goss} for basic facts about Drinfeld modules.

\medskip

\begin{Proofof}{Theorem \ref{main}(b)}

We first reduce to the case where $\End_{\bar{K}}(\phi)=A$. Replacing $K$ by a finite extension if necessary, we may assume that $\End_{\bar{K}}(\phi) = \End_K(\phi)$. Let $R=\End_K(\phi)$, then since $\phi$ has generic characteristic, $R$ is an order in a purely imaginary extension $k'/k$, i.e. $k'$ has only one place above $\infty$. Furthermore, $[k':k]$ divides $r$. Denote by $A'$ the integral closure of $A$ in $k'$. 

By \cite[Prop 4.7.19]{Goss} there exists a Drinfeld $A$-module $\psi$ and an isogeny $P:\phi\to\psi$, defined over $K$, such that $\End_K(\psi) = A'$. Now $P$ induces a morphism $\phi_\tor(L)\to\psi_\tor(L)$, and the dual isogeny $\hat{P}$ likewise induces a morphism $\psi_\tor(L)\to\phi_\tor(L)$, of degree independent of $L$. Hence
\[
c_1|\psi_\tor(L)| \leq |\phi_\tor(L)| \leq c_2|\psi_\tor(L)|
\]
for constants $c_1,c_2>0$ independent of $L$.

Now $\psi$ may be extended to a Drinfeld $A'$-module of rank $r'= r/[A':A]$, which we denote by $\psi'$. We claim that $\psi_\tor(L) = \psi'_\tor(L)$. Let $c\in A'$, then $\psi'_c \in\End_K(\psi)$, and $\hat{\psi}'_c\circ\psi'_c = \psi_d$ for some $d\in A$, where $\hat{\psi}'_c$ denotes the dual of $\psi'_c$ as an isogeny. Hence $\ker(\psi'_c)\subset \ker(\psi_d)$ and so $\psi'_\tor(L)\subset\psi_\tor(L)$. The other inclusion is obvious.

Thus it suffices to prove Theorem \ref{main}(b) with $(\phi,A,r,\gamma)$ replaced by $(\psi',A',r',1/r')$. 
The result now follows from Theorem \ref{GaloisModules} together with the following important result of Pink and R\"utsche \cite{PR}:

\begin{Thm}[Pink-R\"utsche]\label{Pink}
Let $\phi$ be a rank $r$ Drinfeld $A$-module in generic characteristic, defined over the finitely generated field $K$. Suppose that $\End_{K}(\phi)=\End_{\bar{K}}(\phi)=A$. Then there exists a constant $C_0$ depending on $\phi$ and on $K$, such that 
such that the index $I(\Fa)$ of the image of the Galois representation on $\phi[\Fa]\cong (A/\Fa)^r$ is bounded by $C_0$ for all $\Fa\subset A$.
\end{Thm}

\end{Proofof}

\subsection{Elliptic curves}

Let $K$ be a finitely generated field of characteristic zero, and $E/K$ an elliptic curve. For a field $L/K$ we denote by $\End_L(E)$ the ring of endomorphisms of $E$ defined over $L$.

\medskip

\begin{Proofof}{Theorem \ref{main}(a)}

Suppose that $E$ does not have complex multiplication. Then $E$ is a rank 2 $(G_K,\BZ)$-module, and the result follows from Theorem \ref{GaloisModules} with $A=\BZ$ and $r=2$ once we have established

\begin{Thm}[Serre]\label{Serre}
Suppose $E/K$ is an elliptic curve without complex multiplication, defined over a finitely generated field $K$ of characteristic zero. Then there exists a constant $C_0$, depending on $E$ and on $K$, such that the index $I(n)$ of the image of the Galois representation on $E[n]\cong (\BZ/n\BZ)^2$ is bounded by $C_0$ for all $n\in\BZ$.
\end{Thm}

\begin{Proof}
Let $j$ denote the $j$-invariant of $E$ and let $K_1=\BQ(j)$. Then by \cite[\S III, Prop 1.4]{Silverman} there exists an elliptic curve $E'/K_1$ with $j(E')=j$ which becomes isomorphic to $E$ over a finite extension of $K$. It suffices to prove Theorem \ref{Serre} with $E$ replaced by $E'$.

The result holds for $K_1$, in the sense that the cokernel of $\rho_n:\Gal(K_1^\sep/K_1)\to\GL_2(\BZ/n\BZ)$ is bounded independently of $n$. Indeed,
If $K_1$ is a number field then this is Serre's celebrated Open Image Theorem \cite{Serre}, whereas if $j$ is transcendental then the result follows by an older result of Weber \cite[p68]{LangEF}. 

Now let $K_1\subset K_2 \subset K$ such that $K_2/K_1$ is purely transcendental and $K/K_2$ is finite. The result also holds for $K_2$ since $\Gal(K_2^\sep/K_2)\cong\Gal(K_1^\sep/K_1)$. As $\Gal(\Ksep/K)$ is a subgroup of index $[K:K_2]$ in $\Gal(K_2^\sep/K_2)$, it follows that the result also holds for $K$.
\end{Proof}

Next, suppose that $E$ has complex multiplication by an order in the quadratic imaginary field $k/\BQ$. After replacing $K$ by a finite extension and $E$ by a $K$-isogenous elliptic curve if necessary, we may assume that $\End_{\bar{K}}(E) = \End_K(E) = A$ is the maximal order in $k$. Now $E$ is a rank 1 $(G_K,A)$-module, and for any finite extension $L/K$ the torsion points in $E(L)$ with respect to the $A$-module structure coincide with the usual torsion points. Again, the result follows from Theorem \ref{GaloisModules} with $r=1$ together with the following consequence of CM theory.

\end{Proofof}

\begin{Thm}[CM]
Suppose $E/K$ is an elliptic curve defined over a finitely generated field $K$ of characteristic zero. Suppose that $\End_{K}(E) = \End_{\bar{K}}(E) = A$ is the maximal order in a quadratic imaginary field $k$.
Then there exists a constant $C_0$, depending on $E$ and on $K$, such that the index $I(\Fa)$ of the image of the Galois representation on $E[\Fa]\cong A/\Fa$ is bounded by $C_0$ for all $\Fa\subset A$.
\end{Thm}

\begin{Proof}
As before, let $j$ denote the $j$-invariant of $E$, and set $K_1=\BQ(j)$. Since $E$ has complex multiplication, $K_1$ is a number field and the result holds for $K_1$ by \cite[\S4.5]{Serre}. The result now extends to $K$ as in the proof of Theorem \ref{Serre}.
%
\end{Proof}



\section{Discussion}

When $M$ is an elliptic curve with complex multiplication, then the lower bound in Theorem \ref{main} improves the bound $|M_\tor(L)| \geq C [L:K]\sqrt{\log\log[L:K]}$ (for suitable fields $L$) often encountered in the literature (e.g. \cite{CX,HS}).


\subsection{Other applications of Galois modules}

Applying Theorem \ref{GaloisModules} to the rank 1 $(G_\BQ,\BZ)$-module $M=\Gm$, we get the well known result that the number of roots of unity in a number field $L/\BQ$ is bounded by $C[L:\BQ]\log\log[L:\BQ]$, for an absolute constant $C>0$.

One may
also bound the orders of $\Gal(\Ksep/L)$-stable submodules (equivalently, degrees of $L$-rational isogenies) of elliptic curves or Drinfeld modules. 

\begin{Prop}
Let $M$ be a $(G_K,A)$-module of rank $r\geq 2$. Let $L/K$ be a finite extension and $H\subset M$ a $\Gal(\Ksep/L)$-stable cyclic submodule of order $\Fa\subset A$, i.e. $H\cong A/\Fa$. Then 
\[
|H| = |\Fa| \leq C[L:K]^{1/(r-1)},
\] 
where the constant $C$ depends on $M,K,r$ and the index $I(\Fa)$, but not on $L$.
\end{Prop} 

\begin{Proof}
One shows, as in Proposition \ref{Kx}, that the stabilizer of $H$ in $\Aut(M[\Fa])$ has order
\begin{eqnarray*}
|\Stab_{\Aut(M[\Fa])}(H)| & = & |(A/\Fa)^\times|\cdot|\Fa|^{r-1} \cdot|\GL_{r-1}(A/\Fa)| \\
& = & |\Fa|^{r^2-r+1}\prod_{\Fp|\Fa}\left(1-\frac{1}{|\Fp|}\right)^2\left(1-\frac{1}{|\Fp|^2}\right)
\cdots \left(1-\frac{1}{|\Fp|^{r-1}}\right)
\end{eqnarray*}
and hence
\begin{eqnarray*}
[K(H):K] & \geq & \frac{1}{I(\Fa)}\frac{|\GL_r(A/\Fa)|}{|\Stab_{\Aut(M[\Fa])}(H)|} \\
& \geq & \frac{\zeta_A(r)}{I(\Fa)} |\Fa|^{r-1}.
\end{eqnarray*}
The result follows.
\end{Proof}

\begin{Cor}
Suppose $M$ is an elliptic curve without complex multiplication defined over a finitely generated field $K$ of characteristic zero, or that $M$ is a Drinfeld $A$-module of rank $r\geq 2$ in generic characteristic with $\End_{\bar{K}}(M)=A$, defined over the finitely generated field $K$. Then there exists a constant $C>0$, depending on $M,K$ and $r$ such that, for any finite extension $L/K$, the degree of any $L$-rational cyclic isogeny $M\to M'$ is bounded by $C[L:K]^{1/(r-1)}$, (where $r=2$ if $M$ is an elliptic curve).\qed 
\end{Cor}

\subsection{Uniform bounds}

One may ask if the constants in Theorem \ref{main} may be chosen independently of $M$ (this would follow from a uniform bound on $I(\Fa)$). 

When $M$ is a Drinfeld module of rank $1$ this was shown by Poonen \cite[Theorem 8]{Poonen}. For Drinfeld modules of higher rank the existence of an upper bound on $|M_\tor(L)|$ depending only on $r$, $A$ and $[L:K]$ is conjectured by Poonen [{\em loc. cit.}], though there are various partial results, typically with the upper bound depending on primes of bad reduction, see for example \cite{Ghioca, Poonen, Schweizer1}.

When $E$ is an elliptic curve over a number field $K$, uniform upper bounds on $|E_\tor(L)|$ do exist, as shown by Mazur, Kamienny and Merel \cite{ Kamienny, Mazur, Merel}, but these bounds are not yet known to be polynomial in the degree $[L:K]$ in general. When $E$ has everywhere good reduction, then we have the explicit bound $|E_\tor(L)| \leq 1977408[L:\BQ]\log[L:\BQ]$, due to Hindry and Silverman, \cite{HS}. 

If $E$ is an elliptic curve with complex multiplication defined over a number field, one may translate Poonen's proof of \cite[Theorem 8]{Poonen} from rank 1 Drinfeld modules to the rank 1 $(G_K,\End(E))$-module $E$, and one obtains
\begin{equation}
|E_\tor(L)| \leq C[L:\BQ]\log\log[L:\BQ],
\end{equation}
where the constant $C$ depends only on the endomorphism ring $\End(E)$. On the other hand, it follows from \cite{PY,Silverberg} that the exponent of the group $E_\tor(L)$ is bounded by $C[L:\BQ]\log\log[L:\BQ]$ for an absolute constant $C$.

\subsection{Abelian varieties}

Suppose $M$ is an abelian variety of dimension $g$ defined over a number field $K$. Masser has shown that $|M_\tor(L)|\leq C\big([L:K]\log[L:K]\big)^g$, see \cite{Masser1, Masser2}. The exponent $g$ is not optimal in general.

The key to our approach is the independence of fields generated by linearly independent torsion points (Proposition \ref{Indep}), which holds because
\begin{equation}\label{eq:dims}
r\cdot \codim \Fix_G(x) = \dim G,
\end{equation}
when $G=\GL_r$. 

For abelian varieties, the image of Galois is contained in the Mumford-Tate group, which is an algebraic subgroup of $\GSp_{2g}$.
This suggests developing a theory of ``symplectic'' Galois modules. However, as $\dim\GSp_{2g} = 2g^2+g+1$ does not factorize, an identity of the form (\ref{eq:dims}) is not possible. This means that one must explicitly estimate the order of $\Fix_{G(A/\Fa)}(H)$ for submodules $H\subset M_\tor$, which requires more effort. This is done by Hindry and Ratazzi \cite{HR1,HR2}, allowing them to obtain the optimal exponent $\gamma$ for which $|M_\tor(L)| \leq C_{\epsilon}[L:K]^{\gamma+\epsilon}$ holds in various cases (here the constant $C_{\epsilon}$ depends on $\epsilon>0$).

\paragraph{Acknowledgements.} 
The seed for this article was sown in discussions with Marc Hindry and Amilcar Pacheco at the XX Escola de \'Algebra in Rio de Janeiro.


\end{document}